\documentclass{article}%
\usepackage{amssymb}
\usepackage{amsmath}
\usepackage{amsfonts}
\usepackage{graphicx}%
\providecommand{\U}[1]{\protect\rule{.1in}{.1in}}
\newtheorem{theorem}{Theorem}

\newtheorem{condition}[theorem]{Condition}

\newtheorem{corollary}[theorem]{Corollary}

\newtheorem{definition}[theorem]{Definition}

\newtheorem{lemma}[theorem]{Lemma}
\newtheorem{notation}[theorem]{Notation}

\newtheorem{proposition}[theorem]{Proposition}
\newtheorem{remark}[theorem]{Remark}

\newenvironment{proof}[1][Proof]{\noindent\textbf{#1.} }{\ \rule{0.5em}{0.5em}}
\begin{document}

\title{Large time behavior in random multiplicative processes. }
\author{Gy\"{o}rgy Steinbrecher\\Association EURATOM-MECT, \\University of Craiova, Physics Faculty, \\A. I. Cuza Street 13, Craiova-200585,\\Romania, Phone: +40 251 415 077.
\and Xavier Garbet\\CEA, IRFM, F-13108 Saint Paul Lez Durance, France.
\and Boris Weyssow\\Association EURATOM-\'{E}tat Belge,\\Physique Statistique et Plasmas, Universit\'{e} Libre de Bruxelles,\\Campus Plain CP 231, B-1050 Bruxelles, Belgium; EFDA-CSU.}

\begin{abstract}
In a general class of one dimensional random differential equation the
convergence of the distribution function of the solution to stationary state
distribution is studied. In particular it is proved the boundedness
respectively the divergence of the fractional order moments of the solution
below respectively above some critical exponent. This exponent is computed. In
particular models it is the heavy tail exponent. When the equation is linear
this exponent determines a new family of weak topologies (stronger compared to
the classical one), related to the convergence to the stationary state.



\end{abstract}
\maketitle

\section{Introduction}


Discrete or continuous time, one dimensional affine stochastic evolution
equations (ASEE) are studied both in mathematical literature (e.g.
Refs.\cite{Goldie, DeSaportaCnt, DeSaportaDiscr, ShawSG, ShawTimeDep}), as
well as in physical literature in Refs.\cite{Coloids, TakayasuSato, SatoCADP,
SGBW}. This interest in ASEE comes partly from the occurrence of heavy tail
(HT) in the stationary probability distribution functions (PDF), in the models
of physical, economical and biological processes. The ASEE are also related to
reduced models of the self-organized criticality phenomenology in plasma
physics (Ref.\cite{SGBW}), and to the renewal processes (Refs.\cite{Goldie,
DeSaportaCnt, DeSaportaDiscr}).

The connection of the our approach to the linear ASEE models with HT in the
stationary PDF is the following. Denote by $X_{t}(\omega)$ a solution of ASEE
\ and $p_{t}(x):=prob_{\omega}(\left\vert X_{t}(\omega)\right\vert \geq x)$.
Suppose also that the ASEE has a stationary PDF, so we denote
$p(x)=\underset{t\rightarrow\infty}{\lim}p_{t}(x)$. We define the heavy tail
exponent $\beta_{c}$ by the asymptotic estimate $p(x)=O(x^{-\beta_{c}}l(x))$,
where $l(x)$ is a slowly varying function.

\ The occurrence of the HT in the stationary PDF is related to the dynamical
effect (Refs.\cite{ShawTimeDep, SGBW}) that will be studied in this article.
When the stationary PDF of the solution $X_{t}(\omega)$ of the ASEE has HT
with exponent $\beta_{c}$, then for $t\rightarrow\infty$ the fractional order
moments $\mathbb{E}[|X_{t}|^{p}]$ remains bounded for $0<p<\beta_{c}$,
respectively diverges for $p>\beta_{c}$. This is related to the "variance
explosion" phenomena studied recently in the mathematical finance
(Ref.\cite{ShawTimeDep}). In this article we obtain a simple, explicit formula
for $\beta_{c}$.

\ Previously explicit algebraic methods for computing $\beta_{c}$ were
elaborated in some special cases: in the framework of the discrete time models
in Refs.\ \cite{Goldie,TakayasuSato}, with i.i.d. additive and multiplicative
noise, and in the case of i.i.d. additive noise and multiplicative noise
modelled by a finite state Markov process in Ref.\cite{DeSaportaDiscr}.

\ In the continuos time case, with the multiplicative noise modelled by a
finite state Markov process, rigorous foundation of the computation of
$\beta_{c}$ was obtained in Ref.\cite{DeSaportaCnt}. In Ref.\cite{SGBW} the
multiplicative and additive noise were modelled by a superposition of
Ornstein-Uhlenbeck processes. An explicit formula for $\beta_{c}$ was derived
by asymptotic methods. In all of the cases the exponent $\beta_{c}$ is
independent of the additive term.

The ASEE model equation that is considered here is a class of one dimensional
random differential equation (RDE ), which extends previous results from
Ref.\cite{SGBW}, by using new topological vector space methods. The additive
and the multiplicative random terms in our model are stationary processes. The
multiplicative term is a generalization of the stationary Gaussian processes,
having a very general, possibly algebraic, correlation decay. Our results and
those from Ref. \cite{DeSaportaCnt} are complementary.

The main part of the new results are related to the convergence to stationary
state in a class of linear RDE. Nevertheless, the formula for $\beta_{c}$ was
extended to a class of non linear RDE.

Define the subspace of real valued continuos functions $C_{\gamma}%
(\mathbb{R})\subset C(\mathbb{R})$ by the condition $f\in C_{\gamma
}(\mathbb{R})\Longleftrightarrow$ $\left\vert f(x)\right\vert =o(x^{-\gamma}%
)$. For the class of \emph{linear models }considered here, we prove the
existence of the suitable defined weak limit of the PDF of $X_{t}$, in the
sense that there exists a random variable $X_{\infty}$, independent of the
initial conditions, such that if $X_{t}$ is the solution of the RDE and
$0<\gamma<\beta_{c}$, then $\underset{t\rightarrow\infty}{\lim}\mathbb{E}%
\left[  f(X_{t})\right]  =\mathbb{E}\left[  f(X_{\infty})\right]  $ for all
$f\in C_{\gamma}(\mathbb{R})$. The topology on the space of probability
measures, related to this new class of weak limit, is stronger compared to the
classical weak topology.

In particular we have $f(x)=|x+z|^{p}\in C_{\gamma}(\mathbb{R})$ for all
$z\in\mathbb{C}$ when $0<p<\gamma<\beta_{c}$, so for $t\rightarrow\infty$, we
have $\mathbb{E}\left[  |X_{t}+z|^{p}\right]  \rightarrow\mathbb{E}\left[
|X_{\infty}+z|^{p}\right]  $. Moreover we prove that if $p>\beta_{c}$ then
$\mathbb{E}\left[  |X_{t}+z|^{p}\right]  \rightarrow\infty$, for all but
eventually one initial conditions. The expectation values, $\mathbb{E}\left[
|X_{t}+z|^{p}\right]  ,~\mathbb{E}\left[  |X_{\infty}+z|^{p}\right]  $ are
related to a class of generalized $L^{p}$ spaces that includes also the
$p\in\,]0,1[$ case. The proofs uses the geometry and topology of these general
$L^{p}$ spaces. In particular we prove that $X_{\infty}\in L^{p}$ for all
$0<p<\beta_{c}$. When heavy tail exists this property determines the heavy
tail exponent $\beta_{c}$.

We apply these $L^{p}$ methods to the study of nonlinear RDE , when the
nonlinear term has a weak nonlinearity. In this nonlinear model we derive a
weaker results. There exists also a critical exponent $\beta_{c}$, such that
for $0<p<\beta_{c}$ the moments $\mathbb{E}\left[  |X_{t}+z|^{p}\right]  $ are
bounded for large $t$ , for all initial conditions. For $p>\beta_{c}$ and for
sufficiently large initial conditions the moments $\mathbb{E}\left[
|X_{t}+z|^{p}\right]  $ diverges exponentially for large time. The same
formula for $\beta_{c}$ give the demarcation line between the order of the
divergent, respectively bounded fractional moments.

\ In all of the cases $\beta_{c}$ is given by a simple analytic formula in
term of two physically significant parameters of the multiplicative term only,
suggesting some non commutative Central Limit Theorem on the affine group.

The structure of this article is the following. Section \ref{themodel}
provides the notations, the linear model and the main theorem. Section
\ref{SectProofMainTheorm} gives a sequence of propositions leading to the
proof of the main results. The study of the non linear model, by using results
from Section \ref{SectProofMainTheorm} is presented in Section
\ref{sectApplication}. A summary of the results is given in Section
\ref{SectConclusions}.

\section{Description of the linear model and the results.\label{themodel}}



\subsection{Notations and definitions.}


The stochastic processes are defined in a fixed probability space
$\{\Omega,\mathcal{F},P\}$ with expectation value $\mathbb{E}_{\omega
}[f(\omega)]=\mathbb{E}[f]=\int_{\Omega}f(\omega)dP(\omega)$. By $\omega$ will
be denoted a generic element of $\Omega$. Two driving $\mathcal{F}-$measurable
stochastic processes $\{\zeta_{t}(\omega),~\phi_{t}(\omega):\mathbb{R}%
\times\Omega\rightarrow\mathbb{R}^{2}$, and a constant $a>0$ (the instability
threshold) defines formally the linear RDE :
\begin{equation}
\frac{dX_{t}(\omega)}{dt}=-(a+\zeta_{t}(\omega))\,X_{t}(\omega)+\phi
_{t}(\omega)\, \label{9}%
\end{equation}

whose solution can be constructed explicitly (see below). The notations
$X_{t}(\omega)$ or$~X_{t}$, for the solutions of RDE will be reserved. Without
loss of generality we impose $\mathbb{E}[\zeta_{t}]=0$ and deterministic
initial conditions $X_{0}(\omega)\equiv x_{0}$. The argument $\omega$ will be
omitted when no confusion arises. We shall denote by $\Xi$ the set of
solutions of Equation (\ref{9}).

\vskip0.5truecm

\noindent\ We will preserve notation $L^{p}$ also for some unusual Lebesgue
spaces with $0<p<1$ (see Ref.\cite{LpFrechet}, page 75):

\begin{notation}
\label{Not_Lp} \ Let $p>0$ and denote $\sigma_{p}:=\min(1,p)$. We define
$\left\Vert f(\omega)\right\Vert _{p}:=\left(  \mathbb{E}[|f|^{p}]\right)
^{\sigma_{p}/p}$ and $L^{p}:=L^{p}(\Omega,\mathcal{F},P)=\{f\,|\left\Vert
f\right\Vert _{p}<\infty\}$. More explicitly for $p\geq1$ we have the usual
norm $\left\Vert f(\omega)\right\Vert _{p}:=\left(  \mathbb{E}[|f|^{p}%
]\right)  ^{1/p}$, while for $0<p\leq1$ we have the distance to the origin
given by $\left\Vert f(\omega)\right\Vert _{p}:=\mathbb{E}[|f|^{p}]$.
\end{notation}


These complete metric vector spaces were already used in the study of
probability distributions having HT (Ref.\cite{Luschgy}).


\begin{remark}
\label{RemL_p} \ Also for $p\in\,(0,1)$ the topology induced by the distance
$\left\Vert f-g\right\Vert _{p}$ and the vector space structures are
compatible. This follows from the general inequality
\begin{equation}
\left\Vert \alpha f(\omega)+g(\omega)\right\Vert _{p}\leq|\alpha|^{\sigma_{p}%
}\,\left\Vert f\right\Vert _{p}+\left\Vert g\right\Vert _{p}
\label{FormulaGenHolder}%
\end{equation}
where $p>0$ and $f(\omega),~g(\omega)\in L^{p}$. All of the $L^{p}$ spaces are
complete (see page 75 from Ref.\cite{LpFrechet}). It is easy to see that for
$p\in\,(0,1)$ the unit ball is not convex and in nonatomic cases their dual is trivial.
\end{remark}

\vskip0.5truecm

When $0<p\leq1$ the Inequality (\ref{FormulaGenHolder}) results from
$\left\vert a+b\right\vert ^{p}\leq\left\vert a\right\vert ^{p}+\left\vert
b\right\vert ^{p}$.

\vskip0.5truecm

The final convergence results will be formulated in the terms of a subspace
$C_{\gamma}(\mathbb{R})$ of continuos functions on $\mathbb{R}$.

\begin{definition}
\label{DefContFunctsubspace} For any $\gamma>0$ we define the subspace
$C_{\gamma}(\mathbb{R})$ of the space of the continuous functions
$C(\mathbb{R})$ by the condition:
\[
f\in C_{\gamma}(\mathbb{R})\Leftrightarrow\underset{\left\vert x\right\vert
\rightarrow\infty}{\lim}\frac{\left\vert f(x)\right\vert }{\left(
1+\left\vert x\right\vert \right)  ^{\gamma}}=0
\]
and the topology on $C_{\gamma}(\mathbb{R})$ by the norm
\begin{equation}
p_{\gamma}(f):=\underset{x\in\mathbb{R}}{\sup}\frac{\left\vert f(x)\right\vert
}{\left(  1+\left\vert x\right\vert \right)  ^{\gamma}}\label{normGama_f}%
\end{equation}

\end{definition}


\ In the proofs the extension of the space of H\"{o}lder-continuous functions
to the whole real line proves useful, in the study of the tail effects with
extreme delocalization, i.e. $\beta_{c}\ll1$:

\begin{definition}
\label{DefH_alfa} If $\alpha\in\,]0,1]$ then $f(x)\in\mathcal{H}_{\alpha
}\Leftrightarrow\{\exists(c_{1}\geq0)$ ~ such that ~ $\forall(x,y\in
\mathbb{R})~[|f(x+y)-f(x)|\leq c_{1}\,|y|^{\alpha}]\}$.
\end{definition}

\vskip0.5truecm

\ We will preserve the same notation, to define a class of functions useful in
the study of the more localized regimes, when the stationary PDF has
$\beta_{c}>1$:

\begin{definition}
\label{DefH_alfa_large}If $\alpha>1$ then $f(x)\in\mathcal{H}_{\alpha
}\Leftrightarrow\{\exists(g_{i}(x)\in\mathcal{H}_{1},~c_{i}\in\mathbb{C}%
\ )~s.t.~f(x)=\sum_{i}c_{i}\left\vert g_{i}(x)\right\vert ^{\alpha}\}$.
\end{definition}

\begin{remark}
\label{RemPolbound}If $f(x)\in\mathcal{H}_{\alpha}$ then $\left\vert
f(x)\right\vert \leq a+b\left\vert x\right\vert ^{\alpha}$, both for
$\alpha\gtreqqless1$, for some constants $a,b$. The exact values of the
constants $c_{i},a,b$ are irrelevant.
\end{remark}

\begin{remark}
\label{RemTestFunctions}It is elementary to check that if $z\in%
\mathbb{C}
$ , $0<\alpha$, then $f(x)=\left\vert z+x\right\vert ^{\alpha}\in
\mathcal{H}_{\alpha}$.
\end{remark}

\vskip0.5truecm \noindent

By using the spaces $\mathcal{H}_{\alpha}$, $C_{\gamma}(\mathbb{R})$ we define
a parametrized families of weak topologies on the set of probabilistic Borel
measures on $\mathbb{R}$. These topologies, with indexed by a fixed $\rho>0$,
are defined as follows:

\begin{definition}
\label{DefK_beta} Denote $\mathcal{K}_{\rho}:=\underset{\alpha\in
]0,\rho\lbrack}{\bigoplus}\mathcal{H}_{\alpha}$ (finite linear combinations).
For $t\rightarrow\infty$, the sequence $\mu_{t}$ of the Borelian measures on
$\mathbb{R}$ converges to the stationary measure $\mu_{\infty}$ in the
topology $\mathcal{T}(\mathcal{K}_{\rho})$ iff $\ \forall f(x)\in
\mathcal{K}_{\rho}$ we have $\underset{t\rightarrow\infty}{\lim}%
\int_{\mathbb{R}}f(x)\,d\mu_{t}(x)=\int_{\mathbb{R}}f(x)\,d\mu_{\infty}(x)$
\end{definition}


In the final formulation of the results, the asymmetry in the definition of
the spaces $\mathcal{H}_{\alpha}$ will be removed. An apparently stronger
class of topologies are given by the following

\begin{definition}
\label{DefTopolTauGamma} For $t\rightarrow\infty$, the sequence $\mu_{t}$ of
the Borelian measures on $\mathbb{R}$ converges to the stationary measure
$\mu_{\infty}$ in the topology $\mathcal{T}(C_{\gamma})$ iff $\ \forall
f(x)\in C_{\gamma}(\mathbb{R})$ we have $\underset{t\rightarrow\infty}{\lim
}\int_{\mathbb{R}}f(x)\,d\mu_{t}(x)=\int_{\mathbb{R}}f(x)\,d\mu_{\infty}(x)$.
\end{definition}


We observe that the $\mathcal{T}(C_{\gamma})$ is strictly stronger than the
weak topology.


\subsection{Specification of the linear model.}


Without loss of generality and for sake of simplicity, we suppose that a
suitable rescaling of the time variable was performed. Consequently, the
multiplicative noise $\zeta_{t}$ obeys the following

\begin{condition}
\label{CondMultiplNoise} The stochastic process $\zeta_{t}$ is centered and
stationary. There exists a unique "diffusion constant" $D>0$ and a set of
positive constants $K_{p}^{\pm}$, all independent of the times $t$ and $s$,
such that for $p>0$, $s\geq1$ and $s-1\leq t\leq s$ we have:
\begin{equation}
K_{p}^{-}\leq\exp\left(  -D\,p^{2}\,s\right)  \mathbb{E}\left[  \exp\left\{
(1-p)\,Y_{s}-Y_{t}\right\}  \right]  \leq K_{p}^{+} \label{163new}%
\end{equation}
where
\begin{equation}
Y_{t}:=\int_{0}^{t}\zeta_{s}\,ds \label{NEW_Y_t}%
\end{equation}

\end{condition}

The exact values of the constants $K_{p}^{-}$ and $K_{p}^{+}$ are irrelevant.
This Condition is satisfied by a large class of Gaussian processes, as shown
in Proposition\textrm{\ \ref{PropMainExpGaussBound}}.

We will prove that the large time asymptotic behavior of the solution of the
RDE (\ref{9}) is determined by the parameter $a$ from Equation (\ref{9}) and
the constant $D$ from Equation (\ref{163new}). So we introduce the following

\begin{notation}
\label{NotConst} : \emph{The main results are expressed in terms of the
critical exponent}
\begin{equation}
\beta_{c}:=a/D \label{criticalexponent}%
\end{equation}
Another important quantity is $\gamma_{p}:=\sigma_{p}\,D\,(p-\beta_{c})$ with
$p>0$, whose sign changes at $p=\beta_{c}$.
\end{notation}

On the additive noise $\phi_{t}(\omega)$ we impose the following

\begin{condition}
\label{CondAditNoise}

\begin{enumerate}
\item The process $\phi_{t}(\omega)$ is stationary.

\item \label{CondAdNoiseIndep}The processes $\zeta_{t}(\omega),~\phi
_{t}(\omega)$ are independent.

\item \label{CondAdditFinMom}We have $\mathbb{E}\left[  \left\vert \phi
_{t}(\omega)\right\vert ^{p}\right]  \leq m_{p}<\infty$, $\forall
p\in\,]0,\beta_{1}[$, where $\beta_{1}>\max(1,\beta_{c})$.

\item Reversibility: the equality in distribution $\phi_{t}(\omega
)\overset{d}{=}\phi_{-t}(\omega)$.
\end{enumerate}
\end{condition}


From Conditions \ref{CondAditNoise}, item \ref{CondAdditFinMom} results that
the additive noise is allowed to have a HT with exponent larger than
$\max(1,\beta_{c})$.

\vskip0.5truecm \noindent


\subsection{The results, for the linear model.}


Under the previous notations, Conditions\textrm{\ \ref{CondMultiplNoise}
}and\textrm{\ \ref{CondAditNoise}}, we have the following results, that will
be demonstrated in the next sections. First we have the main Theorem and its Corollary


\begin{theorem}
\label{TheoremFinal} Let $f\in C_{\gamma}(\mathbb{R})$ and $0<\gamma<\beta
_{c}$. There exists $X_{\infty}(\omega)\in\bigcap_{p\in]0,\beta_{c}[}L^{p}$
such that:
\end{theorem}

\begin{enumerate}
\item $f(X_{\infty})\in\ L^{1}$

\item We have $\underset{t\rightarrow\infty}{\lim}\mathbb{E}[f(X_{t}%
)]=\mathbb{E}[f(X_{\infty})]$. In particular if $p\in\,]0,\beta_{c}[$ and
$z\in\mathbb{C}$, then $\underset{t\rightarrow\infty}{\lim}\mathbb{E}%
\left\vert X_{t}+z\right\vert ^{p}=\mathbb{E}\left\vert X_{\infty
}+z\right\vert ^{p}<\infty$.

\item If $p>\beta_{c}$ then, except for at most a special value of the initial
condition $X_{t=0}$, we have $\underset{t\rightarrow\infty}{\lim}%
\mathbb{E}\left\vert X_{t}+z\right\vert ^{p}=\infty$.
\end{enumerate}


This Theorem leads immediately to the following Corollary:

\begin{corollary}
\label{CorollaryFinal}Denoting by $F_{t}(x)$ and $F_{\infty}(x)$ the
cumulative PDF of $X_{t}$ and of $X_{\infty}$, respectively, then:

\begin{enumerate}
\item If $\gamma\in\,]0,\beta_{c}[$ and $f\in C_{\gamma}(\mathbb{R})$, then
$\underset{t\rightarrow\infty}{\lim}\int_{\mathbb{R}}f(x)\,dF_{t}%
(x)=\int_{\mathbb{R}}f(x)\,dF_{\infty}(x)$. In particular, for any complex
constant $z$ we have $\underset{t\rightarrow\infty}{\lim}\int_{\mathbb{R}%
}|x+z|^{\gamma}\,dF_{t}(x)=\int_{\mathbb{R}}|x+z|^{\gamma}\,dF_{\infty
}(x)<\infty$.

\item If $\gamma>\beta_{c}$ then except for at most a special value of the
initial condition, we have $\underset{t\rightarrow\infty}{\lim}\int%
_{\mathbb{R}}|x+z|^{\gamma}\,dF_{t}(x)=+\infty$.
\end{enumerate}
\end{corollary}


The Item 1 of this Corollary states that weak convergence of the measure
generated by $F_{t}(x)$ in the topology $\mathcal{T}(C_{\gamma})$.\ 

The previous theorem results from the following technical Lemma and its Corollary



\begin{lemma}
\label{Theor_p_less_beta} There exists $X_{\infty}(\omega)\in\bigcap
_{p\in]0,\beta_{c}[}L^{p}$ such that:
\end{lemma}

\begin{enumerate}
\item $\{f(x)\in\mathcal{K}_{\beta_{c}}\mathcal{\}}\Rightarrow\{f(X_{\infty
})\in\ L^{1}\}$.

\item $\{\ f(x)\in\mathcal{K}_{\beta_{c}}\}\Rightarrow\{\underset{t\rightarrow
\infty}{\lim}\mathbb{E}[f(X_{t})]=\mathbb{E}[f(X_{\infty})]\}$. In particular
for $p\in\,]0,\beta_{c}[$ and $z\in\mathbb{C}$, we have
$\underset{t\rightarrow\infty}{\lim}\mathbb{E}\left\vert X_{t}+z\right\vert
^{p}=\mathbb{E}\left\vert X_{\infty}+z\right\vert ^{p}<\infty$.

\item If $p>\beta_{c}$ then, except for at most a special value of the initial
condition $X_{t=0}$, we have $\underset{t\rightarrow\infty}{\lim}%
\mathbb{E}\left\vert X_{t}+z\right\vert ^{p}=\infty$.
\end{enumerate}

\vskip0.5truecm

This Lemma leads immediately to the following Corollary:

\begin{corollary}
\label{CorolHeavyTail} Denoting by $F_{t}(x)$ and $F_{\infty}(x)$ the
cumulative PDF of $X_{t}$ and of $X_{\infty}$, respectively, then:

\begin{enumerate}
\item If $p\in\,]0,\beta_{c}[$ and $f(x)\in\mathcal{H}_{p}$, then
$\underset{t\rightarrow\infty}{\lim}\int_{\mathbb{R}}f(x)\,dF_{t}%
(x)=\int_{\mathbb{R}}f(x)\,dF_{\infty}(x)$. In particular for any complex
constant $z$ we have $\underset{t\rightarrow\infty}{\lim}\int_{\mathbb{R}%
}|x+z|^{p}\,dF_{t}(x)=\int_{\mathbb{R}}|x+z|^{p}\,dF_{\infty}(x)<\infty$.

\item If $p>\beta_{c}$ then except for at most an eventual special value of
the initial condition we have $\underset{t\rightarrow\infty}{\lim}%
\int_{\mathbb{R}}|x+z|^{p}\,dF_{t}(x)=+\infty$.
\end{enumerate}
\end{corollary}

\noindent


The Item 1 of this Corollary states that weak convergence of the measure
generated by $F_{t}(x)$ in the topology $\mathcal{T}(\mathcal{K}_{\rho})$.\ 

We will prove that the topologies $\mathcal{T}(\mathcal{K}_{\rho})$ are
equivalent $\mathcal{T}(C_{\gamma})$, when restricted to probabilistic measures.

\begin{remark}
\label{M1} The convergence/divergence of the moments for $p\lessgtr\beta_{c}$
does not imply the existence of the HT. Indeed, if $\phi_{t}(\omega)=0$ and
$p\in\,]0,\beta_{c}[$, then $\underset{t\rightarrow\infty}{\lim}\left\Vert
X_{t}\right\Vert _{p}=\left\Vert X_{\infty}\right\Vert _{p}=0$. See Remark
\ref{RemarkX}.
\end{remark}


\begin{remark}
\label{RemExpConvergence}From the proofs will be clear that for $p<\beta_{c}$,
when $\gamma_{p}<0$ (see Notation \ref{NotConst}) the speed of the convergence
for large time of $\left\Vert X_{t}\right\Vert _{p}$ is $O(\exp(\gamma_{p}%
t))$. If $p<\beta_{c}$ then the speed of divergence is $\left\Vert
X_{t}\right\Vert _{p}=O(\exp(\gamma_{p}t))$.
\end{remark}


\section{Proof of the Theorem \ref{TheoremFinal}.\textbf{\ }%
\label{SectProofMainTheorm}}

For technical reasons we introduce the following additional notations

\begin{notation}
\label{NotStochProc}
\begin{align*}
A_{t}  &  :=\exp\left(  -at-Y_{t}\right) \\
B_{t}  &  :=\int_{0}^{t}\phi_{\tau}A_{t}/A_{\tau}\,d\tau\\
H_{t}  &  :=\int_{0}^{t}\phi_{\tau}A_{\tau}\,d\tau\\
d(t)  &  =\mathbb{E}[Y_{t}^{2}]/2
\end{align*}
Furthermore, all the constants denoted by $K_{n}$ with $n$ integer, are
independent of the time variables $t,s,x,y$.
\end{notation}



\vskip0.5truecm.



\subsection{Some initial results.}



\noindent Because the driving noises are classical functions, then according
to the Notations \ref{NotStochProc} the rigorous, integral form of Equation
(\ref{9}), with the initial condition $x_{0}$, is
\begin{equation}
X_{t}=x_{0}\,A_{t}+B_{t}\,. \label{300}%
\end{equation}

\begin{remark}
\label{RemStrictIneq} From Equation \textrm{(\ref{300})} and Equation
\textrm{(\ref{FormulaGenHolder})} we get:
\begin{equation}
\left\vert |x_{0}|^{\sigma_{p}}\left\Vert A_{t}\right\Vert _{p}-\left\Vert
B_{t}\right\Vert _{p}\right\vert \leq\left\Vert X_{t}\right\Vert _{p}%
\leq|x_{0}|^{\sigma_{p}}\left\Vert A_{t}\right\Vert _{p}+\left\Vert
B_{t}\right\Vert _{p}\,. \label{330}%
\end{equation}

\end{remark}

The following proposition is the key point in the difficult part of the proofs
that needs $L^{p}$ bounds when $0<p<1$. This case include also the study of
the "extreme heavy tail" effects (Ref.\cite{SGBW}), when $\beta_{c}$ is very
small.

\begin{proposition}
\label{PropGenFunctionalIneq} Let $U(\omega)\geq0$ and $V(\omega)\geq0$ be
random variables such that $V,~U\,V\in L^{1}$and $\mathbb{E}\left[  V\right]
>0$. Then for $p\in\,]0,1]$ we have
\begin{equation}
\mathbb{E}\left[  U^{p}\,V\right]  \leq\left(  \mathbb{E}\left[  U\,V\right]
\right)  ^{p}\,\left(  \mathbb{E}\left[  V\right]  \right)  ^{1-p}\,.
\label{161}%
\end{equation}

\end{proposition}

\begin{proof}
This result follows from Jensen inequality applied to the concave function
$x\rightarrow x^{p}$. Define a new expectation value $\mathbb{E}_{1}\left[
f\right]  :=\mathbb{E}\left[  fV\right]  /\mathbb{E}\left[  V\right]  $.

For any concave function $g(x):\mathbb{R\rightarrow}\mathbb{R}$ by Jensen
inequality we have $\mathbb{E}_{1}\left[  g(f)\right]  \leq g\left\{
\mathbb{E}_{1}\left[  f\right]  \right\}  $. Consider now $g(x):=x^{p}$ and
$\ f(\omega)=U(\omega)$. We obtain%
\[
\mathbb{E}\left[  U^{p}\,V\right]  /\mathbb{E}\left[  V\right]  \leq\left(
\mathbb{E}\left[  U\,V\right]  /\mathbb{E}\left[  V\right]  \right)  ^{p}%
\]
that leads to the Inequality (\ref{161}).
\end{proof}

\vskip0.5truecm

The following Proposition results from Condition
\textrm{\ref{CondMultiplNoise}} by setting $s=t$:

\begin{proposition}
\label{BasicGaussianInequality} Under the
Condition\textrm{\ \ref{CondMultiplNoise}}, there exist positive constants
$D$, $K_{p}^{+}$, $K_{p}^{-}$, $K_{1}$, and $K_{2}$ such that $\forall
p\in\,]0,\infty\lbrack$ and $t\geq0$ we have
\begin{align}
K_{p}^{-}\,\exp(p^{2}\,D\,t)  &  \leq\mathbb{E}\left[  \exp(-p\,Y_{t})\right]
\leq K_{p}^{+}\exp(p^{2}\,D\,t)\label{162}\\
K_{2}\,\exp\left(  t\,\gamma_{p}\right)   &  \leq\left\Vert A_{t}\right\Vert
_{p}\leq K_{1}\exp\left(  t\,\gamma_{p}\right)  \label{162new}%
\end{align}

\end{proposition}



\begin{remark}
\label{RemarkX} From Inequality (\ref{162new}) and Notation \ref{NotConst}
results that for $0<p<\beta_{c}$ we have $\left\Vert A_{t}\right\Vert
_{p}\rightarrow0$ when $t\rightarrow\infty$. In the particular case $\phi
_{t}(\omega)=0$, by Equation (\ref{300}) we obtain also $\left\Vert
X_{t}\right\Vert _{p}\rightarrow0$.
\end{remark}


\begin{proposition}
\label{PropMainExpGaussBound} Suppose the process $\zeta_{t}$ is Gaussian,
stationary, centered, with continuous realizations and has the correlation
decay for $t\rightarrow\infty$
\begin{equation}
C(\left\vert t\right\vert )=\mathbb{E}(\zeta_{t+\tau}\,\zeta_{\tau
})=\mathcal{O}(1+\left\vert t\right\vert ^{-2-\varepsilon}). \label{162_21}%
\end{equation}
Then Condition \textrm{\ref{CondMultiplNoise}} is satisfied and the
Inequalities\textrm{\ (\ref{163new}, \ref{162}, \ref{162new})} results.
\end{proposition}

\begin{proof}
Because $\zeta_{t}$ is stationary, we will use a method similar to the
derivation of the Taylor-Green-Kubo formula (Refs.\cite{Balescu, Taylor,
Green, Kubo}). For the sake of completeness we adapt the proof. By using
Equations (\ref{NEW_Y_t}, \ref{162_21}) and the Notations \ref{NotStochProc},
we obtain%
\begin{equation}
2\,d(t)=\mathbb{E}[Y_{t}^{2}]=2\,\int\limits_{0}^{t}(t-x)\,C(x)\,dx
\label{nw1}%
\end{equation}
Combined with Equation (\ref{162_21}) it results that for $t\rightarrow\infty$
we have $d(t)=D\,t+\mathcal{O}(1)$, where%
\begin{equation}
D=\int_{0}^{\infty}\,C(|\tau|)\,d\tau=\int_{0}^{\infty}\mathbb{E}[\zeta_{\tau
}\zeta_{0}]d\tau\label{Powerspectrum}%
\end{equation}
Also from Equations (\ref{162_21}, \ref{nw1}) for $\left\vert x\right\vert <1$
we have%
\begin{equation}
d(t+x)-d(t)=\mathcal{O}(1) \label{nw2}%
\end{equation}
$\ $ As a consequence, there exists a constant $K_{3}$ such that for all
$t\geq0$ we have uniformly in $t$: $(D\,t-K_{3})\leq d(t):=\mathbb{E}%
[Y_{t}^{2}]/2\leq(D\,t+K_{3})$.

Because $Y_{t}$ has stationary increments we obtain $\mathbb{E}[Y_{y}%
\,Y_{t}]=d(t)+d(y)-d(y-t)$. With $Z=(1-p)\,Y_{s}-Y_{t}$ we get $\mathbb{E}%
\left[  Z^{2}/2\right]  =p^{2}\,d(s)+[p\,(d(t)-d(s))+(1-p)\,d(s-t)]$.

From $|s-t|\leq1$ and Equation (\ref{nw2}) it results that the last bracket is
uniformly bounded. Since $Z$ is Gaussian, we obtain the Inequality
(\textrm{\ref{163new}}) and by
Proposition\textrm{\ \ref{BasicGaussianInequality}} the Inequalities
\textrm{(\ref{162}, \ref{162new})}.
\end{proof}

\vskip0.5truecm

Observe that by Equation (\ref{Powerspectrum}) the constant $D$ is related to
the zero frequency component of the correlation function. \noindent We shall
prove the following general

\begin{proposition}
\label{PropMaiSmalpBound} Let $p \in\, ]0, 1[$, $0\leq y-1\leq x\leq y$ and
suppose that

\begin{enumerate}
\item The stochastic processes $Y_{t},\zeta_{t}$ satisfies
Condition\textrm{\ \ref{CondMultiplNoise}.}

\item The stochastic process $f_{t}$ is stationary, independent of $Y_{t}$ and
we have the bound uniformly in $t$:
\begin{equation}
\left\vert \mathsf{E}[f_{t}]\right\vert <K_{4}<\infty\,. \label{1664}%
\end{equation}

\end{enumerate}

Denote
\begin{equation}
b_{f}(x,y):=\left\Vert \int_{x}^{y}f_{t}\,A_{t}\,dt\right\Vert _{p}
\label{1665}%
\end{equation}
Then for some constant $K_{5}$ we have
\begin{equation}
b_{f}(x,y)\leq K_{5}\exp\left(  y\gamma_{p}\right)  \,. \label{difficultbound}%
\end{equation}

\end{proposition}

\begin{proof}
Since for $p\in\,]0,1[$ we have $\sigma_{p}=p$ and $\left\Vert g\right\Vert
_{p}=\mathbb{E}\left[  \left\vert g\right\vert ^{p}\right]  $, the Equation
(\ref{1665}) can be rewritten as $b_{f}(x,y)=\mathbb{E}\left[  U^{p}%
\,V\right]  $, where
\begin{align*}
V(\omega) &  =\exp\left[  -p\,Y_{y}(\omega)\right]  \\
U(\omega) &  =\int_{x}^{y}f_{t}(\omega)\,A_{t}(\omega)\,\exp\left[
Y_{y}(\omega)\right]  \,dt\,.
\end{align*}
Recalling Proposition \ref{PropGenFunctionalIneq} and Inequality (\ref{161})
we get:
\begin{equation}
b_{f}(x,y)\leq\left(  \mathbb{E}\left[  U\,V\right]  \right)  ^{p}\,\left(
\mathbb{E}\left[  V\right]  \right)  ^{1-p}\,.\label{a28}%
\end{equation}
Then, on the one hand, according to Proposition \ref{BasicGaussianInequality}
and Inequality (\ref{162}), the term $\mathbb{E}\left[  V\right]  $ in
Equation (\ref{a28}) is bounded by
\begin{equation}
\mathbb{E}\left[  V\right]  \leq K_{6}\,\exp\left(  p^{2}\,y\,D\right)
\,.\label{a29}%
\end{equation}
On the other hand, due to the independence of $f_{t}$ and $Y_{t}$ we obtain
\[
\mathbb{E}\left[  U\,V\right]  =\int_{x}^{y}\mathbb{E}[f]\,\mathbb{E}%
[A_{t}\,\exp\left\{  (1-p)Y_{y}\right\}  ]\,dt\,.
\]
Recalling Inequality (\ref{1664}), and the definition of $A_{t}$, it follows
that
\[
\left\vert \mathbb{E}\left[  U\,V\right]  \right\vert \leq K_{4}\,\int_{x}%
^{y}\,\exp(-a\,t)\mathbb{E}\left[  \exp\left\{  (1-p)\,Y_{y}-Y_{t}\right\}
\right]  \,dt
\]
which, with the help of Inequality (\ref{163new}), reduces to
\begin{equation}
\left\vert \mathbb{E}\left[  U\,V\right]  \right\vert \leq K_{7}%
\,\exp(-a\,y+p^{2}\,D\,y)\,.\label{bx1}%
\end{equation}
After simple calculations, recalling Inequalities (\ref{a28}, \ref{a29},
\ref{bx1}), the Inequality (\ref{difficultbound}) is obtained.
\end{proof}

\begin{remark}
\label{REMStatIncrements} From the Condition\textrm{\ \ref{CondAditNoise}},
for fixed $T$, by symmetry and stationarity, we have an equality in
distribution of the direct and reflected process $\phi_{t}\overset{d}{=}%
\phi_{T-t}$.
\end{remark}

\vskip0.5truecm

\noindent From the previous Remark, it follows the following

\begin{proposition}
\label{PropBtEqHt} If $f(x)\in\mathcal{H}_{p}$ then $\mathbb{E}[f(B_{t}%
)]=\mathbb{E}[f(H_{t})]$ and in particular $\left\Vert B_{t}+z\right\Vert
_{p}=\left\Vert H_{t}+z\right\Vert _{p}$, where $z\in%
\mathbb{C}
,~p>0$.
\end{proposition}

\begin{proof}
From the Definitions \ref{DefH_alfa} \ref{DefH_alfa_large} results, that if
$f(x)\in\mathcal{H}_{p}$ then $f(x)$ is continuos, hence measurable. From
Remark \ref{RemPolbound} we have\ $\left\vert f(x)\right\vert \leq c_{1}%
+c_{2}\left\vert x\right\vert ^{p}$. According to the Notations
\ref{NotStochProc}
\[
B_{t}=\int_{0}^{t}\phi_{\tau}\exp\left(  -a(t-\tau)-(Y_{t}-Y_{\tau})\right)
\,d\tau
\]
From Condition \ref{CondMultiplNoise} results that the process $Y_{\tau}$ has
stationary increments, so for fixed $t$ we have the equality in distribution
of the processes $Y_{t}-Y_{\tau}\overset{d}{=}Y_{t-\tau}$. On the other hand,
from Condition \ref{CondAdNoiseIndep} results that because the processes
$\phi_{\tau},Y_{\tau}$ are independent, we have the equality in distribution%
\[
\left(  \phi_{\tau},Y_{t}-Y_{\tau}\right)  \overset{d}{=}\left(  \phi_{\tau
},Y_{t-\tau}\right)
\]
Consequently
\[
B_{t}\overset{d}{=}\int_{0}^{t}\phi_{\tau}\exp\left(  -a(t-\tau)-Y_{t-\tau
}\right)  \,d\tau
\]

By the change of the integration variable $\tau\rightarrow t-\tau$ and from
Remark \ref{REMStatIncrements} results the equality in distribution
\begin{equation}
B_{t}\overset{d}{=}\int_{0}^{t}\phi_{t-\tau}\exp\left(  -a\tau-Y_{\tau
}\right)  \,d\tau\label{pr5}%
\end{equation}
According to the Remark \ref{REMStatIncrements}, \emph{for fixed }$t$\emph{
}we obtain $\phi_{\tau}\overset{d}{=}\phi_{t-\tau}$. By the independence of
the processes $\phi_{\tau},\zeta_{\tau}$, consequently the independence of
$\phi_{\tau},Y_{\tau}$ we have%
\[
\left(  \phi_{\tau},Y_{\tau}\right)  \overset{d}{=}\left(  \phi_{t-\tau
},Y_{\tau}\right)
\]
So from Equation (\ref{pr5}) results
\[
B_{t}\overset{d}{=}\int_{0}^{t}\phi_{\tau}\exp\left(  -a\tau-Y_{\tau}\right)
\,d\tau
\]
or (see Notation \ref{NotStochProc}) $B_{t}\overset{d}{=}H_{t}$ which
completes the proof.
\end{proof}

\vskip0.5truecm \noindent The following Lemma and its
Corollary\textrm{\ \ref{Corol_f(c_t)}} will be used in subsequent works.

\begin{lemma}
\label{LemmaCauchyBound} Under the previous Notations \ref{NotConst},
\ref{NotStochProc}, Conditions \textrm{\ref{CondMultiplNoise}} and
\textrm{\ref{CondAditNoise}}, we have

\begin{enumerate}
\item \label{x1l1} If $v\geq u\geq0$ , $\gamma_{p}\neq0$ \ and $0<p<\beta_{1}$
then we have the $L^{p}$ bound
\begin{equation}
\left\Vert H_{v}-H_{u}\right\Vert _{p}\leq K_{8}\,\left\{  \exp(u\,\gamma
_{p})+\exp(v\,\gamma_{p})\right\}  \label{350}%
\end{equation}
for some constant $K_{8}$ independent of $u,v$.

\item \label{x1l2} There exists $X_{\infty} \in\bigcap_{0 < q < \beta_{c}} \,
L^{q} $ such that $\underset{t \rightarrow\infty}{\lim}\left\Vert
H_{t}-X_{\infty}\right\Vert _{p}=0$ when $p \in\, ]0,\beta_{c}[$.

\item \label{x113} $\forall p\in]0,\beta_{1}[$ we have%
\begin{equation}
\left\Vert B_{t}\right\Vert _{p}\leq K_{8}\,\left[  1+\exp(t\,\gamma
_{p})\right]  \label{BoundBt}%
\end{equation}

\end{enumerate}
\end{lemma}

\begin{proof}
\begin{enumerate}
\item \emph{If} $p\in\,[1,\beta_{1}[$ \emph{then} (see Notation
\ref{NotStochProc}) we obtain $\left\Vert H_{v}-H_{u}\right\Vert _{p}\leq$
$\int_{u}^{v}\left\Vert \phi_{s}\,A_{s}\right\Vert _{p}ds$. Recalling
Conditions \ref{CondAditNoise} items 2 and 3, it follows that%
\[
\left\Vert \phi_{s}\,A_{s}\right\Vert _{p}=\left\Vert \phi_{s}\right\Vert
_{p}\left\Vert A_{s}\right\Vert _{p}\leq m_{p}\,\left\Vert A_{s}\right\Vert
_{p}.
\]
Consequently, by using the Inequality (\textrm{\ref{162new}}) results
\[
\left\Vert H_{v}-H_{u}\right\Vert _{p}\leq m_{p}K_{1}\,\int_{u}^{v}\exp\left(
t\,\gamma_{p}\right)  ds
\]
which immediately lead to Equation (\ref{350}).

\ \ \ In the "very heavy tail" case when $p\in]0,1[$ we use the Proposition
\ref{PropMaiSmalpBound}, where we replace $f_{t}\equiv\phi_{t}$ and Equation
(\ref{1665}). We introduce the notations $[x]$ for the integer part of $x$,
the notations $n_{v}:=[v-u],v_{-}:=u+n_{v}$ and use the decomposition
\[
\int_{u}^{v}\phi_{\tau}A_{\tau}d\tau=\left(  \int_{u}^{u+1}+\cdot\cdot
\cdot+\int_{u+n_{v}-1}^{u+n_{v}}+\int_{u+n_{v}}^{v}\right)  \phi_{\tau}%
A_{\tau}d\tau
\]
We get%
\[
\left\Vert H_{v}-H_{u}\right\Vert _{p}=b_{\phi}(u,v)\leq b_{\phi}%
(v_{-},v)+\sum_{k=1}^{n_{v}}b_{\phi}(u+k-1,u+k)
\]
By using Proposition \ref{PropMaiSmalpBound}, Equation (\ref{difficultbound})
with $f_{t}\equiv\phi_{t}$, after simple estimations we obtain Equation
(\ref{350}), which completes the proof.

\item Let $p\in\,]0,\beta_{c}[$. Then we have $\gamma_{p}<0$ and from Equation
(\ref{350}) it results that if $t_{n}\rightarrow\infty$ then $H_{t_{n}}$ is a
Cauchy sequence in $L^{p}$.

The $L^{p}$ spaces, including $p\in]0,1[$, are complete (Ref.\cite{LpFrechet},
page 75), therefore an $X_{\infty}(\omega)\in L^{p}$ exists such that
$\left\Vert H_{t}-X_{\infty}\right\Vert _{p}\rightarrow0$. Because $0<p\leq
q\Rightarrow L^{p}\supset L^{q}$, it results that $X_{\infty}(\omega
)\in\bigcap_{0<q<\beta_{c}}L^{q}$.

\item It results from Proposition\textrm{\ \ref{PropBtEqHt}} and Inequality
\textrm{(\ref{350})}, by setting $u=0,~v=t$.
\end{enumerate}
\end{proof}

\begin{lemma}
\label{LemmaCont_Lp_L1} Let $f(x)\in\mathcal{H}_{p}$ fixed. Then $\Psi\in
L^{p}\Rightarrow f(\Psi)\in L^{1}$.The corresponding application
$~L^{p}\rightarrow\mathbb{C}$, defined as $L^{p}\ni\Psi\rightarrow
\mathbb{E}\left[  f(\Psi)\right]  \in\mathbb{C}$, is continuous.
\end{lemma}

\begin{proof}
From $f(x)\in\mathcal{H}_{p}$ results $|f(x)|\leq a+b\,|x|^{p}$ (see Remark
\ref{RemPolbound}). Hence $\Psi\in L^{p}\Rightarrow f(\Psi)\in L^{1}$.

For the proof of the continuity we consider first the case $0<p\leq1$. Let
$\Psi,\chi\in L^{p}$ and recall that in this case $\ f(x)\in\mathcal{H}%
_{p}\Rightarrow|f(x+y)-f(x)|\leq c_{1}\,|y|^{p}$. Then it follows that
\[
\left\vert \mathbb{E}\left[  f(\Psi+\chi)\right]  -\mathbb{E}\left[
f(\Psi)\right]  \right\vert \leq\mathbb{E}\left\vert f(\Psi+\chi
)-f(\Psi)\right\vert \leq c_{1}\mathbb{E}\left\vert \chi\right\vert ^{p}%
=c_{1}\left\Vert \chi\right\Vert _{p}%
\]
which proves the continuity in the case of the (unusual) $L^{p}$ space.

In the case $p\geq1$ we observe that according to the Definition
\ref{DefH_alfa_large} it is sufficient to prove the continuity for the set of
functions of the form $f(x)=\left\vert g(x)\right\vert ^{p}$ when$~g(x)\in
\mathcal{H}_{1}$, that generates the space $\mathcal{H}_{p}$ by finite linear combinations.

Let $\Psi,\chi\in L^{p}$ with $p\geq1$ and recall that in this case
$\ f(x)\in\mathcal{H}_{p}\Rightarrow|g(x+y)-g(x)|\leq c_{1}\,|y|$. So it is
sufficient to prove that
\begin{equation}
\left\Vert \chi\right\Vert _{p}\rightarrow0\Rightarrow\mathbb{E}\left[
f(\Psi+\chi)\right]  \rightarrow\mathbb{E}\left[  f(\Psi)\right]  \label{nw3}%
\end{equation}
\ In the case $f(x)=\left\vert g(x)\right\vert ^{p}$ the previous Equation
(\ref{nw3}) is equivalent to the condition:
\begin{equation}
\left\Vert \chi\right\Vert _{p}\rightarrow0~\Rightarrow\left\Vert g(\Psi
+\chi)\right\Vert _{p}\rightarrow\left\Vert g(\Psi)\right\Vert _{p}
\label{nw4}%
\end{equation}
Observe that from the usual ($p\geq1$) H\"{o}lder inequality and from
$g(x)\in\mathcal{H}_{1}$ we obtain
\[
\left\vert \left\Vert g(\Psi+\chi)\right\Vert _{p}-\left\Vert g(\Psi
)\right\Vert _{p}\right\vert \leq\left\Vert g(\Psi+\chi)-g(\Psi)\right\Vert
_{p}\leq\left\Vert c_{1}\chi\right\Vert _{p}%
\]
which completes the proof for the case $p\geq1$.
\end{proof}

\begin{corollary}
\label{Corol_f(c_t)} Let $p\in\,]0,\beta_{c}[$ and $f(x)\in\mathcal{H}_{p}$.
Then $\mathbb{E}[f(H_{t})]\rightarrow\mathbb{E}[f(X_{\infty})$ and
$\mathbb{E}[f(B_{t})]\rightarrow\mathbb{E}[f(X_{\infty})]$ for $t\rightarrow
\infty$.
\end{corollary}

\begin{proof}
The convergence of $\mathbb{E}[f(H_{t})$ results from Lemma
\ref{LemmaCauchyBound} part \ref{x1l2}, and Lemma \ref{LemmaCont_Lp_L1}. From
Proposition \ref{PropBtEqHt} results $\mathbb{E}[f(B_{t})]=\mathbb{E}%
[f(H_{t})]$, which completes the proof.
\end{proof}

\begin{remark}
\label{RemHT} The distribution of $X_{\infty}\in L^{p}$, with $p\in
\,]0,\beta_{c}[$ is identical with $\lim_{t\rightarrow\infty}H_{t}:=\int%
_{0}^{\infty}\phi_{\tau}\,A_{\tau}d\tau$, where the limit is in $L^{p}$.
Clearly, in generic cases $X_{\infty}$ is non degenerate, and for the limiting
measure $\mu(x)=prob(\left\vert X_{\infty}(\omega)\right\vert \geq x)$, we
obtain $0<\int_{0}^{\infty}y^{p}d\mu\left(  y\right)  <\infty$ for all
$p\in]0,\beta_{c}[$.
\end{remark}


\subsection{Proof of\ the Lemma \ref{Theor_p_less_beta}.}


\begin{proof}
The random variable $X_{\infty}\in\cap_{0<p<\beta_{c}}L^{p}$ was
\textbf{constructed} in Lemma \ref{LemmaCauchyBound} part \ref{x1l2}. Let
$p\in\,]0,\beta_{c}[$ and $f(x)\in\mathcal{H}_{p}$. If $f(x)\in\mathcal{K}%
_{\beta_{c}}$, then it can be represented in the form $f(x)=\sum_{i=1}%
^{n}f_{i}(x)$ where $f_{i}(x)\in\mathcal{H}_{p_{i}}$ for $p_{i}\in
\,]0,\beta_{c}[$ and $1\leq i\leq n$. It is therefore sufficient to consider
the case when $f(x)\in\mathcal{H}_{p}$ for $p\in\,]0,\beta_{c}[$.

\begin{enumerate}
\item It is sufficient to prove that $\{f(x)\in\mathcal{H}_{p},p\in
\,]0,\beta_{c}[\}\Rightarrow\{f(X_{\infty})\in L^{1}\}$. We now use
$X_{\infty}\in L^{p}$ and Lemma \ref{LemmaCont_Lp_L1}.

\item Let $t\rightarrow\infty$. It is sufficient to prove that $\{f(x)\in
\mathcal{H}_{p},p\in\,]0,\beta_{c}[\}\Rightarrow\{\underset{t\rightarrow
\infty}{\lim}\mathbb{E}[f(X_{t})]=\mathbb{E}[f(X_{\infty})]\}$. From
Proposition \ref{PropBtEqHt} it results that%
\[
\mathbb{E}[f(B_{t})]=\mathbb{E}[f(H_{t})].
\]
Consequently, according to Corollary \ref{Corol_f(c_t)}, we have
\[
\underset{t\rightarrow\infty}{\lim}\mathbb{E}[f(B_{t})]=\mathbb{E}%
[f(X_{\infty})]<\infty\,.
\]
We will use now Lemma \ref{LemmaCont_Lp_L1} with $\chi=x_{0}\,A_{t}$ and
$\Psi=B_{t}$. Because for $0<p<\beta_{c}$ we have $\exp(t\gamma_{p}%
)\rightarrow0$ then, by Inequality (\textrm{\ref{162new}}) results $\left\Vert
A_{t}\right\Vert _{p}\rightarrow0$, so $\left\Vert \chi\right\Vert
_{p}\rightarrow0$. Consequently by Lemma \ref{LemmaCont_Lp_L1} we obtain%
\[
\mathbb{E}[f(X_{t})]=\mathbb{E}[f(B_{t}+x_{0}\,A_{t})]\rightarrow
\mathbb{E}[f(X_{\infty}).
\]
$\mathbb{E}[f(X_{t})]=\mathbb{E}[f(B_{t}+x_{0}\,A_{t})]\rightarrow
\mathbb{E}[f(X_{\infty})$. The last remark follows from the fact that
$f(x)\in\mathcal{H}_{p}$ for all $z\in\mathbb{C}$ if the function $f(x)$ is
defined by $f(x):=\left\vert x+z\right\vert ^{p}$.

\item Consider the case when $p>\beta_{c}$ \textit{i.e.} $\gamma_{p}>0$ and
therefore $\exp(\gamma_{p}t)\rightarrow+\infty$. Suppose, ad absurdum, that
there exist two initial conditions $x_{0}$ and $x_{0}^{\prime}$ such that the
corresponding solutions $X_{t}(\omega)$ and $X_{t}^{\prime}(\omega)$ of
Equation (\ref{9}) are bounded in $L^{p}$.

The stochastic process $Z_{t}(\omega)=X_{t}(\omega)-X_{t}^{\prime}(\omega)$ is
also bounded and it is determined by the homogenous equation. According to
Equation (\ref{300}), it results that%
\[
Z_{t}(\omega)=(x_{0}-x_{0}^{\prime})\,A_{t}(\omega)
\]
and thus $\left\Vert Z_{t}\right\Vert _{p}=\left\vert (x_{0}-x_{0}^{\prime
})\right\vert ^{\sigma_{p}}\left\Vert A_{t}\right\Vert _{p}$. Then, recalling
the first part of Inequality\textrm{\ (\ref{162new})} we finally obtain%
\begin{equation}
\left\Vert Z_{t}\right\Vert _{p}\geq\left\vert (x_{0}-x_{0}^{\prime
})\right\vert ^{\sigma_{p}}\,K_{2}\,\exp(\gamma_{p}t) \label{111}%
\end{equation}
Because when $p>\beta_{c}$ \textit{we have.} $\gamma_{p}>0$ and therefore
$\exp(\gamma_{p}t)\rightarrow+\infty$., by Equation (\ref{111}) the proof is completed.
\end{enumerate}
\end{proof}


\subsection{Proof of the Theorem \ref{TheoremFinal}.}

The main point of the proof is the density of the space $\mathcal{H}_{\alpha}$
in the space $C_{\gamma}(\mathbb{R})$. For the convenience of the reader, \ we
recall some very general definitions and a generalization by Nachbin of the
Kakutani-Stone density theorem from (Ref. \cite{Nachbin}, Theorem 2 ) to the
non compact case. See also (Ref. \cite{Kashimoto}).

\subsubsection{The Nachbin approximation Lemma}

We use the terminology and notations from (Ref. \cite{Kashimoto}). Let $\ X$ a
completely regular topological space. By $C(X)$ we denote the space of
continuous real valued functions on $X$.

A directed set $V$, where $V\subset C(X)$, is a set such that for any
$v_{1},v_{2}\in V$ there exist $\lambda>0$ and $v\in V$ such that
$v_{1}(x)\leq\lambda v(x)$ and $v_{2}(x)\leq\lambda v(x)$.

We suppose also that for any $x\in X$ there exists $v\in V$ such that $v(x)>0$
(i.e. the set $V$ is pointwise strictly positive). A function $f(.)\in C(X)$
is an element of the weighted space $CV_{\infty}(X)$ if and only if for all
$v\in V$ the function $v(x)f(x)$ vanishes at infinity.

The topology on $CV_{\infty}(X)$ is defined by the seminorms indexed by
elements of $V$, denoted $p_{v}(f)$%
\[
p_{v}(f)=\underset{x\in X}{\sup}\left\vert v(x)f(x)\right\vert
\]

where $f\in CV_{\infty}(X),~v\in V$.

\ Recall, a lattice $L\subset C(X)$ is a set of functions closed under $\min$
and $\max$
\[
f,g\in L\Rightarrow\min\left[  f(x),g(x)\right]  \in L~and~\max\left[
f(x),g(x)\right]  \in L~
\]
The following Lemma of Nachbin (Refs \cite{Nachbin,Kashimoto}) will be used


\begin{lemma}
\label{NachbinLemma} (Nachbin \cite{Nachbin}) Let $X$ be a completely regular
space, $V$ a pointwise strictly positive set of weights, $L$ a sublattice of
$CV_{\infty}(X)$ and $f\in CV_{\infty}(X)$. Then $\ f$ can be approximated by
elements of $L$ in the $CV_{\infty}(X)$ topology if and only if for any
$x,y\in X$ and $\varepsilon>0$ there exists $g\in L$ such that%
\begin{align}
\left\vert f(x)-g(x)\right\vert  &  <\varepsilon\label{Nachbin1}\\
\left\vert f(y)-g(y)\right\vert  &  <\varepsilon\label{Nachbin2}%
\end{align}

\end{lemma}


The following Corollary results

\begin{corollary}
\label{CorollaryDensityC_gamma} For any $0<\alpha<\gamma$ the space
$\mathcal{H}_{\alpha}$ is dense in $C_{\gamma}(\mathbb{R})$ in the topology
induced by the norm from Equation (\ref{normGama_f}).
\end{corollary}

\begin{proof}
In order to apply the Lemma \ref{NachbinLemma} in the our case, we consider
$X=\mathbb{R}$ and \ the set $V$ consists of a single function $v(x)=\left(
1+\left\vert x\right\vert \right)  ^{-\gamma}$. In this case the space
$CV_{\infty}(X)$ is $C_{\gamma}(\mathbb{R})$ from Definition
\ref{DefContFunctsubspace}. First we recall Remark \ref{RemPolbound} and
observe that $\mathcal{H}_{\alpha}$ is contained in $C_{\gamma}(\mathbb{R})$
when $0<\alpha<\gamma$.

To apply Lemma \ref{NachbinLemma}, for any fixed $\alpha>0$ we define the
subset $L$ of $\mathcal{H}_{\alpha}$ as follows: $h\in L$ if and only if there
exists $\eta\in\mathcal{H}_{1}$ such that for any $x\in\mathbb{R}$ we have
$h(x)=\left\vert \eta(x)\right\vert ^{\alpha}$. It is easy to see that in this
case $h\in\mathcal{H}_{\alpha}$ for any $\alpha>0$.

Indeed, for $\alpha\geq1$ it follows from the very definition of
$\mathcal{H}_{\alpha}$. For $0<\alpha\leq1$ we have $\left\vert \left\vert
\eta(x)\right\vert ^{\alpha}-\left\vert \eta(y)\right\vert ^{\alpha
}\right\vert \leq\left\vert \eta(x)-\eta(y)\right\vert ^{\alpha}$ and now we
use that \ $\eta\in\mathcal{H}_{1}$, i.e. $\left\vert \eta(x)-\eta
(y)\right\vert \leq c\left\vert x-y\right\vert $. So we obtain $\left\vert
h(x)-h(y)\right\vert =\left\vert \left\vert \eta(x)\right\vert ^{\alpha
}-\left\vert \eta(y)\right\vert ^{\alpha}\right\vert \leq\lbrack c\left\vert
x-y\right\vert ]^{\alpha}$

It is easy to verify that the set $L$ defined below is a lattice.

Indeed, suppose that $h_{1},h_{2}\in L\subset\mathcal{H}_{\alpha}$, for some
fixed $\alpha>0$. This means that they can be represented as $h_{1}%
(x)=\left\vert \eta_{1}(x)\right\vert ^{\alpha},~h_{2}(x)=\left\vert \eta
_{2}(x)\right\vert ^{\alpha}$, with $\eta_{1,2}\in\mathcal{H}_{1}$. Then we
have $h(x)=\min\left[  h_{1}(x),h_{2}(x)\right]  =\left\vert \eta
(x)\right\vert ^{\alpha}$ , where $\eta(x)=\min[\eta_{1}(x),\eta_{2}%
(x)]\in\mathcal{H}_{1}$, because $\mathcal{H}_{1}$ is closed under $\min
,~\max$ operations. In a similar manner we prove that $L$ is closed under
$\max[]$.

In order to use the Lemma \ref{NachbinLemma} it remains to prove the
restrictions given by Equations (\ref{Nachbin1}, \ref{Nachbin2}).

Moreover, because the positive and negative part of a function from
$C_{\gamma}(\mathbb{R})$ belongs to $C_{\gamma}(\mathbb{R})$ too, it is
sufficient to verify that any non negative function $f(x)\geq0$ from
$C_{\gamma}(\mathbb{R})$ can be approximated by elements of $L$.

But it is easy to check that by using the family of functions $g(x)=\left\vert
ax+b\right\vert ^{\alpha}\in L$, $a,b\in\mathbb{C}$ , the conditions given by
Equations (\ref{Nachbin1}, \ref{Nachbin2}) of the previous Lemma
\ref{NachbinLemma} are obeyed, for $f(x)\geq0$, which completes the proof.
\end{proof}


\subsubsection{Proof of the convergence of the measures (Theorem
\ref{TheoremFinal})}

\begin{proof}
We will use Lemma \ref{Theor_p_less_beta} and its Corollary
\ref{CorolHeavyTail}. In order to complete the proof, by extending the results
from Lemma \ref{Theor_p_less_beta} and Corollary \ref{CorolHeavyTail} to
Theorem \ref{TheoremFinal}, it is sufficient to prove that form the
convergence of the measures in the weak topology $\mathcal{T}(\mathcal{K}%
_{\beta_{c}})$ it follows the convergence of the measures in all of the
topologies $\mathcal{T}(C_{\gamma})$ where $\gamma<\beta_{c}$ .

In this end it is sufficient to prove that from $\underset{t\rightarrow
\infty}{\lim}\int_{\mathbb{R}}f(x)\,dF_{t}(x)=\int_{\mathbb{R}}%
f(x)\,dF_{\infty}(x)$ for $f(x)\in\mathcal{H}_{\alpha}$, $\alpha\in
\,]0,\beta_{c}[$, where $F_{t}(x),F_{\infty}(x)$ where defined in Corollary
\ref{CorolHeavyTail}, it follows that for any $g(.)\in C_{\gamma}%
(\mathbb{R}),~0<\alpha<\gamma<\beta_{c}$ we have the same convergence%
\[
\underset{t\rightarrow\infty}{\lim}\int_{\mathbb{R}}g(x)\,dF_{t}%
(x)=\int_{\mathbb{R}}g(x)\,dF_{\infty}(x)
\]
Denote $\nu_{t}\left(  x\right)  :=F_{t}(x)-F_{\infty}(x)$ . From Corollary
\ref{CorolHeavyTail} results that there exists some finite constant
$K_{\gamma}$ such that%
\begin{equation}
\int_{\mathbb{R}}\left(  1+\left\vert x\right\vert \right)  ^{\gamma
}\left\vert d\nu_{t}(x)\right\vert \,=K_{\gamma},~0<\gamma<\beta_{c}
\label{convmes1}%
\end{equation}

\emph{We emphasize that this is the point where }$\gamma<\beta_{c}$\emph{ is
used. } Let $g(.)\in C_{\gamma}(\mathbb{R}),$ and according to the previous
Corollary \ref{CorollaryDensityC_gamma} consider the sequence $f_{n}%
(x)\in\mathcal{H}_{\alpha}$ of approximants $g(x)$ in the topology
of$~C_{\gamma}(\mathbb{R})$ where $0<\alpha<\gamma<\beta_{c}$. So we have%
\begin{equation}
\underset{n\rightarrow\infty}{\lim}\underset{x\in\mathbf{R}}{\sup}\left\vert
\frac{g(x)-f_{n}(x)}{\left(  1+\left\vert x\right\vert \right)  ^{\gamma}%
}\right\vert =\underset{n\rightarrow\infty}{\lim}p_{\gamma}(g-f_{n})=0
\label{convmes2}%
\end{equation}
We have to prove that for any $\varepsilon>0$ there exists $T_{\varepsilon}$
such that it $t\geq T_{\varepsilon}$ then $\left\vert \int_{\mathbb{R}%
}g(x)d\nu_{t}\left(  x\right)  \right\vert \leq\varepsilon$. We have%
\[
\int_{\mathbb{R}}g(x)d\nu_{t}\left(  x\right)  =\int_{\mathbf{R}}f_{n}%
(x)d\nu_{t}\left(  x\right)  \,+\int_{\mathbf{R}}\frac{g(x)-f_{n}(x)}{\left(
1+\left\vert x\right\vert \right)  ^{\gamma}}\left(  1+\left\vert x\right\vert
\right)  ^{\gamma}d\nu_{t}\left(  x\right)
\]
or%
\[
\left\vert \int_{\mathbb{R}}g(x)d\nu_{t}\left(  x\right)  \right\vert
\leq\left\vert \int_{\mathbf{R}}f_{n}(x)d\nu_{t}\left(  x\right)
\,\right\vert +p_{\gamma}(g-f_{n})\int_{\mathbf{R}}\left(  1+\left\vert
x\right\vert \right)  ^{\gamma}\left\vert d\nu_{t}\left(  x\right)
\right\vert
\]

By Equations (\ref{convmes1}, \ref{convmes2}) we select $n$ such that
$p_{\gamma}(g-f_{n})\leq\varepsilon/(2K_{\gamma})$.

Because $f_{n}(x)\in\mathcal{H}_{\alpha}$, $\alpha\in\,]0,\beta_{c}[$ from
Corollary \ref{CorolHeavyTail} results that there exists $T_{\varepsilon}$
such that if $t\geq T_{\varepsilon}$ then $\left\vert \int_{\mathbb{R}}%
f_{n}(x)d\nu_{t}\left(  x\right)  \right\vert \leq\varepsilon/2$ .

This last inequality completes the proof of the convergence in the topology
$\mathcal{T}(C_{\gamma})$ , $\gamma<\beta_{c}$, under the hypothesis that we
have convergence in the topology $\mathcal{T}(\mathcal{K}_{\beta_{c}})$.

Consequently the results from\ Lemma \ref{Theor_p_less_beta} and its Corollary
\ref{CorolHeavyTail} can be extended to Theorem \ref{TheoremFinal} and
Corollary \ref{CorollaryFinal} which completes the proof.
\end{proof}


\section{Application \label{sectApplication}}

We consider now a generalization of the Equation (\ref{9}), containing a new
non linear term $\Psi_{t}\left[  X_{t}(\omega),\omega\right]  $. We write the
equation symbolically as%
\begin{equation}
\frac{dX_{t}(\omega)}{dt}=-(a+\zeta_{t}(\omega))\,X_{t}(\omega)+\Psi
_{t}\left[  X_{t}(\omega),\omega\right]  \label{apl1}%
\end{equation}

The driving noise terms $\zeta_{t}(\omega),~\Psi_{t}\left[  X_{t}%
(\omega),\omega\right]  $ satisfy the following restrictions.

\begin{condition}
\label{CondNlinTermApl}The stochastic process $\zeta_{t}(\omega)$ obeys
Condition \ref{CondMultiplNoise}. There exists a stochastic process $\phi
_{t}(\omega)\geq0$ such that
\begin{equation}
\left\vert \Psi_{t}\left[  X_{t}(\omega),\omega\right]  \right\vert \leq
\phi_{t}(\omega)\text{ a.e.} \label{apl2}%
\end{equation}

and $\phi_{t}(\omega)$ obeys Condition \ref{CondAditNoise}.
\end{condition}


In this case the existence and uniqueness of the solution is guaranteed by
Theorem $221$, page $65$ from (Ref. \cite{LArnold}). The heavy tail effect is
manifest now by the boundedness/divergence of the moments $\mathbb{E}\left[
|X_{t}+z|^{p}\right]  \ $when $p\lessgtr\beta_{c}$:

\begin{theorem}
Denote by $X_{t}(\omega)$ the solution of the Equation (\ref{apl1}). Under
Condition \ref{CondNlinTermApl} we have the following behavior for
$t\rightarrow\infty$.
\end{theorem}

\begin{enumerate}
\item when $0<p<\beta_{c}$ the fractional moment $\mathbb{E}\left[
|X_{t}+z|^{p}\right]  $ are uniformly bounded in $t$.

\item For $p>\beta_{c}$ and for sufficiently large values of $X_{0}$, the
moments $\mathbb{E}\left[  |X_{t}+z|^{p}\right]  $ diverges.
\end{enumerate}

\begin{proof}
In analogy to the Equation (\ref{300}) and Notations \ref{NotStochProc} we
have the rigorous, implicit, integral form of the Equation (\ref{apl1})%
\begin{align}
X_{t}  &  =x_{0}\,A_{t}+\widetilde{B_{t}}\label{apl3}\\
\widetilde{B_{t}}  &  :=\int_{0}^{t}\Psi_{t}\left[  X_{t}(\omega
),\omega\right]  A_{t}/A_{\tau}\,d\tau\label{apl4}%
\end{align}

Because $\left\Vert f(\omega)\right\Vert _{p}:=\left(  \mathbb{E}%
[|f|^{p}]\right)  ^{\sigma_{p}/p}$ (see Notation \ref{Not_Lp}) it is
sufficient to study the\emph{ boundedness} of $\left\Vert X_{t}+z\right\Vert
_{p}$, and by Equation (\ref{FormulaGenHolder}), \emph{the boundedness of
}$\left\Vert X_{t}\right\Vert _{p}$\emph{.}

From Equations (\ref{apl3}, \ref{FormulaGenHolder}) results%
\begin{equation}
\left\vert x_{0}\right\vert ^{\sigma_{p}}\left\Vert \,A_{t}\right\Vert
_{p}-\left\Vert \,\widetilde{B_{t}}\right\Vert _{p}\leq\left\Vert
\,X_{t}\right\Vert _{p}\leq\left\vert x_{0}\right\vert ^{\sigma_{p}}\left\Vert
\,A_{t}\right\Vert _{p}+\left\Vert \,\widetilde{B_{t}}\right\Vert _{p}
\label{apl5}%
\end{equation}

According to the Inequalities (\ref{162new}, \ref{apl5}), we have%
\begin{equation}
\left\vert x_{0}\right\vert ^{\sigma_{p}}K_{2}\,\exp\left(  t\,\gamma
_{p}\right)  -\left\Vert \,\widetilde{B_{t}}\right\Vert _{p}\leq\left\Vert
\,X_{t}\right\Vert _{p}\leq\left\vert x_{0}\right\vert ^{\sigma_{p}}K_{1}%
\exp\left(  t\,\gamma_{p}\right)  +\left\Vert \,\widetilde{B_{t}}\right\Vert
_{p} \label{apl5a}%
\end{equation}
So it is sufficient to prove (in analogy to Inequality (\ref{BoundBt})) the
new bound%
\begin{equation}
\left\Vert \,\widetilde{B_{t}}\right\Vert _{p}\leq K_{8}\,\left[
1+\exp(t\,\gamma_{p})\right]  \label{NewBoundBt}%
\end{equation}

Observe that because $A_{t}/A_{\tau}>0$, from Equations (\ref{apl2},
\ref{apl4}) we have the inequality
\begin{equation}
\left\vert \widetilde{B_{t}}(\omega)\right\vert \leq\int_{0}^{t}\phi
_{t}(\omega)A_{t}/A_{\tau}\,d\tau\label{apl6}%
\end{equation}

According to the Condition \ref{CondNlinTermApl}, in the Equation (\ref{apl6})
we can identify%
\[
\int_{0}^{t}\phi_{t}(\omega)A_{t}/A_{\tau}\,d\tau=B_{t}%
\]
.Finally, by using Inequality (\ref{BoundBt}) we obtain the Inequality
(\ref{NewBoundBt}). From \ Inequalities (\ref{apl5a}, \ref{NewBoundBt} ), it
results%
\begin{equation}
\left\Vert \,X_{t}\right\Vert _{p}\leq\left\vert x_{0}\right\vert ^{\sigma
_{p}}K_{1}\exp\left(  t\,\gamma_{p}\right)  +K_{8}\,\left[  1+\exp
(t\,\gamma_{p})\right]  \label{apl7}%
\end{equation}
respectively%
\begin{equation}
\left\vert x_{0}\right\vert ^{\sigma_{p}}K_{2}\,\exp\left(  t\,\gamma
_{p}\right)  -K_{8}\,\left[  1+\exp(t\,\gamma_{p})\right]  \leq\left\Vert
\,X_{t}\right\Vert _{p} \label{apl8a}%
\end{equation}

If $p<\beta_{c}$ then we have $\gamma_{p}<0$\ and from Inequality (\ref{apl7})
results that $\left\Vert \,X_{t}\right\Vert _{p}$ is bounded. If $p>\beta_{c}$
and if the initial conditions are sufficiently large, such that $\left\vert
x_{0}\right\vert ^{\sigma_{p}}>K_{8}/K_{2}$, then we have $\gamma_{p}>0$ so by
Inequality (\ref{apl8a}) results that $\left\Vert \,X_{t}\right\Vert _{p}$
diverges exponentially, which completes the proof.
\end{proof}



\section{Conclusions.\label{SectConclusions}}


In a class of one dimensional random differential equations the large time
behavior of the solution was studied. When the equation is linear we proved
that the convergence to stationary state can be described in the framework of
new class of weak topologies on the set of probability distribution of the
solution. These topologies are stronger, compared to the classical weak
topology. The study of the weak convergence involves in a natural way the
study of \emph{the convergence of the fractional order moments of the
solution}. A critical exponent $\beta_{c}$ was defined such that the moments
of the solution, of order $p$ remains bounded if $p<\beta_{c}$ and diverges,
on a massive set of initial conditions, when $p<\beta_{c}$. When heavy tail
exists then $\beta_{c}$ is the heavy tail exponent. The speed of
convergence/divergence, for large time, of the moments of order $p$ of the
solution is exponential (see Remark \ref{RemExpConvergence}), depending on
$\ p-\beta_{c}$

\ The strength of the new family of topologies, that describe the approach to
the steady state, increases with $\beta_{c}$.

\ An important result is the exact and simple Equation (\ref{criticalexponent}%
) for $\beta_{c}$, in term of the parameters from the multiplicative term
only. The convergence in this topology of the distribution function to the
stationary distribution was proved.

A new topological vector space method was used in the proofs. By these new
methods we obtained an exact formula on the critical exponent $\beta_{c}$ also
in the case of a class of nonlinear models described by Equation (\ref{apl1}).
Generalization to higher dimensions remains a challenging open problem.

\section{Acknowledgements}

This work, supported by the European Communities under the contract of
Association between EURATOM and MEdC, EURATOM and CEA, respectively
EURATOM-ULB, was carried out within the framework of the European Fusion
Development Agreement. The views and opinions expressed herein do not
necessarily reflect those of the European Commission. S.G. acknowledges ULB
and CEA-Cadarache for warmth hospitality and to C. P. Niculescu from
University of Craiova, Romania, for helpful discussions.




\end{document}